\documentclass[review,12pt]{article}
\usepackage[utf8]{inputenc}
\usepackage[english]{babel}
\usepackage[T1]{fontenc}
\usepackage{orcidlink,thumbpdf,lmodern}
\usepackage{amsmath}
\usepackage{amsfonts}
\usepackage{framed}
\usepackage{multirow}
\usepackage{float}
\usepackage{hyperref}
\usepackage{authblk}
\usepackage{natbib}

\newtheorem{theorem}{Theorem}

\begin{document}

\title{rSRD: An \textsf{R} package for the Sum of Ranking Differences statistical procedure}

\date{}

\author[1,2]{Balázs R. Sziklai~\orcidlink{0000-0002-0068-8920}\thanks{email: sziklai.balazs@krtk.hun-ren.hu}}
\affil[1]{HUN-REN Centre for Economic and Regional Studies, Budapest, Hungary}
\affil[2]{Department of Operations Research and Actuarial Sciences, Corvinus University of Budapest, Budapest, Hungary}

\author[3]{Attila Gere~\orcidlink{0000-0003-3075-1561}}
\affil[3]{Hungarian University of Agriculture and Life Sciences, Budapest, Hungary}

\author[4]{Károly Héberger~\orcidlink{0000-0003-0965-939X}}
\affil[4]{Plasma Chemistry Reasearch Group, Institute of Materials and Environmental Chemistry,\\ HUN-REN Research Centre for Natural Sciences, Budapest, Hungary}

\author[5]{Jochen Staudacher~\orcidlink{0000-0002-0619-4606}}
\affil[5]{Fakult\"at Informatik, Hochschule Kempten, Kempten, Germany}

\maketitle

\begin{abstract}
Sum of Ranking Differences (SRD) is a relatively novel, non-para-metric statistical procedure that has become increasingly popular recently. SRD compares solutions via a reference by applying a rank transformation on the input and calculating the distance from the reference in $L_1$ norm. Although the computation of the test statistics is simple, validating the results is cumbersome---at least by hand. There are two validation steps involved. Comparison of Ranks with Random Numbers, which is a permutation-test, and cross-validation combined with statistical testing. Both options impose computational difficulties albeit different ones.
The rSRD package was devised to simplify the validation process by reducing both validation steps into single function calls.
In addition, the package provides various useful tools including data preprocessing and plotting.  The package makes SRD accessible to a wide audience as there are currently no other software options with such a comprehensive toolkit. This paper aims to serve as a guide for practitioners by offering a detailed presentation of the features.
\end{abstract}

\section[Introduction]{Introduction} \label{sec:intro}

Sum of Ranking Differences first appeared in the field of analytical chemistry \citep{hebergerSumRankingDifferences2010,KOLLARHUNEK2013139}. In chemometrics, comparing methods or models via a reference is quite natural, as materials are often tested against some industry standard, a benchmark provided by some accredited laboratory, or some known theoretical value.

The idea caught on quickly, as the use of references is also common in other fields. Recommendation systems utilize consumers' preferences to generate content. In choice modeling, the consensus ranking or ground truth is considered the ideal or true ranking of the evaluated items. In web searches, a queried image can serve as a reference. A predefined set of features also plays a role in pattern recognition and computer vision. In feature matching, the reference consists of a set of features or attributes considered important or representative in the context of the problem. Consequently, SRD has a wide range of applications, from machine learning \citep{Moorthy2017} and multi-criteria decision-making \citep{PETCHROMPO2022} to network science \citep{sziklai_lengyel_2024}, outlier detection \citep{brownfield2017consensus}, pharmacology \citep{Vajna2012}, political science \citep{sziklaiApportionmentDistrictingSum2020}, ecology \citep{WATROBSKI20221}, and even sports \citep{West2018,Gere2022}. 

The test statistics of SRD is computed as follows. A rank transformation is performed on the input. Both in the solution vectors and in the reference, values are replaced by ranks. The smallest value of each vector obtains a rank of 1, while the largest $n$ (the size of the input). Ties are resolved by fractional ranking, that is, tied values are replaced by the arithmetic mean of their corresponding ranks. After rank transformation, the distances between the solutions and the reference are computed. SRD uses the $L_1$-norm. We discuss the advantage of this distance metric as well as possible alternatives in Section~\ref{sec:lit_distance}.


SRD values are normalized by the largest possible difference between two rankings of size $n$. In this way, SRD values of different problems can be compared with each other. The obtained values already rank the solutions---the one with the smallest SRD score is the closest to the reference, hence the best. However, the real benefit comes from the validation options.

SRD values are validated in two ways. The first option is the Comparison of Ranks with Random Numbers (CRRN). Informally, it is a permutation-test measuring the probability that a random ranking produces a distance as extreme as the solution's SRD score. By convention, we accept solutions (i.e.\ reject $H_0$) below the 5\% significance threshold. Between 5-95\% solutions might have only artificial resemblance to the reference, but in fact they are not distinguishable from a random ranking. Above 95\% the solutions seem to rank the objects in a reverse order---again $H_0$ is rejected with 5\% significance, although the solution is not accepted since it is closer to the exact opposite of the reference. Thus, SRD can detect significant similarity as well as dissimilarity.

The difficulty of this option is that there is more than one way to assign a distribution to the SRD values. If neither the solution nor the reference ranking contains ties and the number of elements is sufficiently high, the distribution is asymptotically normal. Should any of these criteria (no ties, high $n$) not be met, the distribution can change substantially. Unfortunately, ties are quite common in applications and so is a small sample size. The derivation of thresholds can only be achieved by simulating the empirical distribution, since in such cases, the SRD values do not follow any widely recognized distribution.


The second validation option is cross-validation combined with statistical testing. Although the SRD values rank the solutions, it is unclear what equal or similar SRD scores signify. Just because two solutions are of the same distance from the reference it doesn't mean they are close to each other. During cross-validation we take random samples from the rows and calculate the SRD scores for each sample. Then, with the help of a statistical test (e.g.\ Wilcoxon,  Alpaydin or Dietterich-test) we can determine whether the two rankings induced by the solution vectors come from the same distribution, and if not, which one should be ranked lower in the SRD ranking. The choice of statistical test has a huge impact on type 1 and type 2 errors. Given the trade-off between the two error rates, ultimately only the practitioner can decide which one is best suited for her use case.

Due to the popularity of the procedure there are various software implementations available. Dávid Bajusz published a python based application uploaded to github\footnote{\url{https://github.com/davidbajusz/srdpy}}, see also \citep{GERE2021}. John Kalivas' homepage offers a MATLAB code\footnote{\url{https://www.isu.edu/chem/faculty/staffdirectoryentries/kalivas-john.html}}. Attila Gere implemented a platform independent shiny application\footnote{\url{ https://attilagere.shinyapps.io/srdonline/}}. Finally, Klára Kollárné Hunek and Károly Héberger developed an MS Excel macro\footnote{\url{http://aki.ttk.mta.hu/srd}}. While these existing packages are undoubtedly valuable, they often lack thorough documentation, and none of them provide the comprehensive toolset found in rSRD. Our package excels in various aspects, offering a broad range of features for data preprocessing, SRD computation, validation, and plotting. Furthermore, it outperforms previous implementations in terms of scalability and precision. As a result, the rSRD package emerges as a much-needed solution, bridging the gap for users of this statistical test. In the following, we formally introduce the SRD test, describe the functions featured in the package, and demonstrate their use through examples.


\section{Related literature} \label{sec:Lit}

In this section, we review the literature based on the three components of SRD analysis: rank transformation, reference ranking, and distance measures employed in testing. In addition, we briefly present the theoretical work that has been done related to SRD's validation options.

\subsection{Ranking transformation in statistical testing}

Sum of Ranking Differences is conceptually similar to Kendall’s Tau \citep{Kendall_1938} and Spearman’s Rank Correlation Coefficient \citep{Spearman1904}, both of which measure ordinal associations between variables. However, there are a few notable differences.

Firstly, SRD uses a fixed reference, whereas Kendall's Tau or Spearman's Rho are typically employed in either pairwise or multivariate comparisons. SRD is a proper distance metric, where 0 indicates perfect positive concordance, and 1 represents perfect negative concordance. In contrast, both Kendall's Tau and Spearman's Rho range from $-1$ to 1, where the scales are opposite, with $-1$ signifying the largest possible discrepancy and 1 representing perfect correlation.

The most significant difference lies in how the statistics are computed. In the case of SRD, distance is measured in the $L_1$-norm, while Kendall's Tau uses the number of inversions (also called Kemeny-distance or Kendall Tau distance), and Spearman's Rho computes the product-moment correlation coefficient. In the absence of ties, this coefficient resembles the Euclidean distance between rankings, in the sense that it uses the squares of rank differences in its computation.

Robustness can be enhanced by applying ranking transformation, especially when the data may not meet the assumptions of parametric methods. Ranking helps mitigate the impact of outliers and skewed distributions, making the analysis more robust. For a comprehensive mathematical discussion the reader is referred to \citep{kendall1970rank}. 

\subsection{Reference values in statistical testing}


The use of reference comes natural in many applications. In econometrics there is often a base period or control group to which data is compared. In recommendation systems, user preferences \citep{Bhowmik2017}, while in web searches, the queried subject \citep{Dwork2001} can serve as a reference. Sometimes a reference is not given externally but is drawn from data. This happens, for instance, in marketing and advertisement research where consumer preferences are aggregated to create the ideal ranking \citep{Lin2010,Haghani2021,Schuster2024}. In text summarization, the goal is to generate a condensed representation of a document while preserving the significant content \citep{Joshi2022}. Web metasearches also apply this technique \citep{Akritidis2022}. The most popular statistical models are the Mallows \citep{mallowsNONNULLRANKINGMODELS1957}, the Thurstone \citep{thurstone1927law} and Placket-Luce \citep{luce1959individual,Plackett1975}.

The Thurstone model was developed in experimental psychology to handle cases when subjects have to arrange a series of stimuli in absolute rank order according to the sensation it prompted. The distribution of sensations from a particular stimulus is assumed to be normal. Furthermore, it assumes equal standard deviations of sensations corresponding to stimuli and equal correlations between pairs of stimuli sensations \citep{Mosteller1951}.

The Plackett-Luce model is a probabilistic model used to represent the likelihood of observed rankings. Plackett was originally motivated by a problem of calculating the odds of racehorses, in particular the probability of whether a horse finishes in the top 3 \cite{Plackett1975}. Holman and Marley proved that if the underlying random variables in Thurstone's approach have an extreme value distribution, the resulting choice probabilities are given by the Plackett-Luce model as well \citep{Diaconis1988}.

Finally, in the Mallows model, there is a central ranking $\pi_0$ and a scale parameter $\lambda$. The highest probability is assigned to $\pi_0$, with probabilities diminishing geometrically as the distance from $\pi_0$ increases. A higher value of $\lambda$ leads to a distribution that is increasingly peaked around $\pi_0$. Mallows initially focused on two specific distance metrics: Spearman's rho distance and Kendall's tau distance.

Mallows, Plackett-Luce, and Thurstone are burdened by both assumptions about the distribution of the underlying data and computational difficulties. In the Thurstone model, non-linear least squares approach is used to estimate the means of the normal distributions, while the standard deviations are assumed to be unit \citep{Lin2010}. In the Mallows model finding the Maximum Likelihood Estimate of $\pi_0$ is the Kemeny’s consensus ranking problem, known to be NP-hard \citep{tang2019mallows}.

Fortunately, there are excellent software packages available to address these issues.  The \emph{PlackettLuce} package provides functions for preparing ranking data to fit the Plackett-Luce model or Plackett-Luce trees \citep{PlackettLuce_vgn}. \emph{BayesMallows} offers a Bayesian version of the Mallows rank model \citep{BayesMallows_vgn}. The \emph{rankdist} package serves as a general platform for distance-based ranking models, including Mallows' \citep{Qian2019}.

Finally, it's important to note that the shared feature among these ranking models is the assumption of a reference ranking. However, the objectives of the Mallows, Plackett-Luce, and Thurstone models differ considerably from those of the Sum of Ranking Differences. Consequently, they cannot serve as substitutes for SRD and \emph{vice versa}.

\subsection{Distance metrics in ranking}\label{sec:lit_distance}

A commonly used distance metric is Kendall tau, quantifying the number of inversions between two rankings. In simpler terms, it measures the number of adjacent transpositions required to transform one ranking into another. The attractiveness of this metric lies in its inherent connection to rankings; it does not make sense with real valued arrays. Both the Mallows model and, unsurprisingly, the Kendall Tau Rank Correlation test uses the Kendall Tau measure.

The $L_1$-norm, sometimes also referred to as Manhattan-distance, was also proposed early on. Following Spearman's work, the $L_1$-norm of rankings, in the absence of tied values, is commonly referred to as Spearman's footrule. A comprehensive exploration of the statistical properties of the footrule is provided by ref.\ \cite{Diaconis_Graham_1977}. Empirical evidence suggests that $L_1$-norm is a sensible choice, comparable to other commonly used distance metrics. Héberger and Škrbić found that SRD (which uses the $L_1$-norm) is slightly stricter than Spearman rho and Kendall tau, that is, it rejects more models \cite{HEBERGER2012}. Sipos et al. on the other hand compared SRD, Kendall tau, Caley distance, and a combination of Caley and SRD, and found that SRD and Kendall tau gives virtually the same results \cite{Sipos2018}.

In various applications, such as sports or information retrieval, the top of a ranking often carries more significance than the bottom. Arguably, the accuracy of the first page of search results is far more critical for a search engine than the accuracy of subsequent pages. Similarly, sports enthusiasts tend to prioritize the top three positions over the last three. Consequently, distance metrics considering the top-$k$ positions have been developed for rank aggregation purposes \citep{Fagin2003}. Another solution that takes position into consideration without discarding some parts of the ranking, is weighting. Weighted distance-metrics for rankings are considered for instance in refs.\ \citep{Gere2022,CHATTERJEE201847} and \citep{sziklai2022efficacy}. However, imposing a weighting always adds a subjective element to the analysis as different weightings will benefit different solutions. Moreover, the weighting will deteriorate the power of statistical tests. Thus, weighted distances should be adopted with care. Finally, incomplete rankings and their aggregation are also studied \citep{RODRIGO2024111882}.

\subsection{SRD valildation}\label{sec:lit_distrib}

Comparison of Ranks with Random Numbers (CRRN) is a type of permutation test. The latter is a test where the $p$-value is the proportion of data configurations yielding a test statistic as extreme as the value observed in the research results \citep{edgington2007randomization}. In our case, this translates to the proportion of permutations that are at most as far from the reference as the solution's ranking. 

If there are no ties in the ranks, the $L_1$-norm is commonly referred to as  Spearman's footrule. Initially dismissed by Kendall \cite{kendall1970rank} due to a lack of known statistical properties, later research by Diaconis and Graham \cite{Diaconis_Graham_1977} proved that the distribution composed of footrule distances is asymptotically normal. In combinatorics, this distance is known as the total displacement of a permutation, and its distribution is of independent interest. Generating the exact distribution is computationally intractable, as the problem is equivalent to counting weighted Motzkin paths of a given area \citep{Guay_Paquet_2014}. The distribution for $n \leq 50$ is available in the On-Line Encyclopedia of Integer Sequences \citep{oeis2024}. However, if ties are present SRD scores do not follow any special distribution and its density curve becomes jagged.

CRRN enables the grouping of solutions, distinguishing acceptable methods similar to the reference from those that are not. The second validation option in SRD is cross-validation, aiding in identifying the true order of solutions with the same or very close SRD scores. Originally, \citet{KOLLARHUNEK2013139} proposed the Wilcoxon signed-rank test for cross-validation purposes. The rSRD package implements two additional statistical tests: the Dietterich ${\mathit{t}}$-test~\citep{dietterichApproximateStatisticalTests1998} and Alpaydin's $\mathcal{F}$-test~\citep{alpaydinCombinedCvTest1999}, both popular cross-validation tools in machine learning. \cite{Sziklai2024} compared the performance of all three tests on synthetic and real data under various parameterizations and input sizes. The best-performing method was the Wilcoxon test with 8 folds, which is the default method of the package. Nevertheless, rSRD provides all three tests as they drastically differ in type I and type II error rates. Although the Wilcoxon-test proved to be a bit
too sensitive in type I scenarios, it is the only method that performed well in
type II situations. On real data, the advantage of Alpaydin and Dietterich methods in type I cases diminishes.


\section{SRD} \label{sec:SRD}

\subsection{Getting and installing rSRD} \label{sec:SRD_math}

The rSRD package is downloadable from the Comprehensive R Archive Network (CRAN) at \url{https://CRAN.R-project.org/package=rSRD} with this paper referring to package version 0.1.8. For efficiency, computationally intensive tasks like cross-validation or calculating the SRD distribution that corresponds to given data are implemented in \textsf{C++}. The \textsf{C++} code is interfaced to R using the established package Rccp by Dirk Eddelbuettel \citep{eddelbuettel2013seamless}. This setup implies that Windows users need to have Rtools with a version greater than or equal to 4.2 installed to build rSRD from source. Furthermore, rSRD imports functionality from the R packages dplyr, ggplot2, ggrepel, gplots, janitor, rlang, stringr and tibble.


\subsection{Mathematical description of SRD} \label{sec:SRD_math}

The input of the SRD-test is an $n\times m$ matrix $A=[a_{ij}]$. The first $m-1$ columns represent the models or methods (\emph{solutions}) that we would like to compare, while rows represent measurements or features (\emph{objects}). The last column has a special role, it contains the reference values for each row. This can be an external reference: a gold standard, a benchmark value, or a previous measurement. In the absence of a known gold standard, the reference can be extracted from the data. This step is referred sometimes to as data fusion or preference aggregation depending on the field. A common solution is to create a new column by taking the average of the row values for each row. The underlying idea is that the random errors of the measurements follow normal distribution and cancel each other out. In the presence of outliers, the row median is also a sensible choice. Furthermore, if the row medians form a ranking then it is the ranking that minimizes the total distance from the solution rankings in $L_1$-norm \citep{Dwork2002}, in other words, it is the ranking closest to the solutions. Depending on the use cases, the row minimum or maximum can also serve as a reference.

After the reference is fixed, a ranking transformation is performed on the input matrix. Values in each column are replaced by ranks. The smallest value receives a rank of 1, the second smallest gets a rank of 2, and so forth. The largest value is replaced by $n$. Ties are resolved by fractional ranking, tied values are replaced by the arithmetic mean of their corresponding ranks. For instance, if $m_{ij}=m_{kj}$ are tied for the 6th and 7th places, they are each replaced by a rank of 6.5. Let us denote the rank-transformed matrix by $R=[r_{ij}]$. The SRD score of the $j$th solution is the distance between the $j$th column and the reference column in $L_1$-norm

$$SRD_j=\sum_{i=1}^n |r_{ij} - r_{im}|.$$

SRD values are normalized by the maximum distance between rankings of size $n$, which has the following explicit formula

\begin{gather}f(n) = \left\lfloor\frac{n^2}{2}\right\rfloor =
\begin{cases}
\frac{n^2}{2} & \mbox{if } n \mbox{ is even} \label{eq:normalization}\\
\frac{n^2-1}{2} & \mbox{if } n \mbox{ is odd}
\end{cases}\end{gather}

Normalized SRD, denoted by $\underline{SRD}_i= SRD_i/f(n)$ is a number between 0 and 1, where 0 means that the solution produces the same ranking as the reference, while 1 means the solution ranks the objects in reverse order.

\subsection{Validation options}

SRD scores already establish a ranking between the solutions. With the validation steps, we are able to extract more information.

One immediate question is whether there are unsuitable solutions. In many applications, the objective is not solely to identify the single best option but rather to have the flexibility to select from a set of solutions that meet certain criteria or requirements. In other scenarios, solutions may represent possible underlying factors of the reference ranking, and our goal is to identify the relevant ones. In such cases, we are interested in distinguishing the good solutions from the bad ones—those that have no significant connection to the reference. Comparison of Ranks with Random Numbers (CRRN) aims to address this issue.

Another interesting question is how reliable the SRD ranking is or how to group the solutions. If there is a huge gap between the SRD scores of the two solutions, by all likelihood their order is correct. But what happens if the SRD values fall close to each other? In particular, what should we do with tied values? To determine whether two rankings are essentially the same or not cross-validation is applied combined with statistical testing (CVST).

\subsubsection{CRRN}

One validation option of SRD is Comparison of Ranks with Random Numbers (CRRN), where SRD scores of the solutions are compared with those of random rankings. The permutation test, a standard technique in statistics, helps determine whether the results occurred by chance. By analogy, consider two binary sequences, $A = 0000000000$ and $B = 1010100011$. While binary sequences of the same length have an equal probability of emerging as a result of a random coin toss, we might still intuitively feel that sequence $A$ is less random than $B$. From a statistical standpoint, our intuition is not entirely unfounded. Coin tosses with very few heads are rare, and the longer the sequence, the rarer they become. For instance, the probability of a ten-bit long sequence containing less than 3 heads is around 0.055, slightly higher than the usual significance threshold. In this example, a person would rightfully be cautious about accepting sequence $A$ as a result of fair coin tosses.


The concept behind CRRN is that a solution closely related to the reference will rank the objects more or less the same way. The less association they have, the less likely they are to rank objects similarly. To illustrate, consider a preference elicitation setting where users rate products or services (e.g. movies). Two users may have similar tastes but different evaluation habits; for instance, one may be stricter than the other. Despite having a large absolute difference in item scores, they might still rank items in the exact same way. Discordant tastes will result in different rankings. If the tastes have no relation to each other, we expect the corresponding ranking to show no pattern. Therefore, the distance between the solution and the reference ranking should be comparable to that of a random ranking and the reference ranking. This can be formulated as a null-hypothesis:\\

$H_0$: \emph{The distance between a solution ranking and the reference ranking is not significantly different from the distance between a random ranking and the reference ranking.}\\

A valuable feature of SRD is that it measures both similarity and dissimilarity. Consequently, we can reject $H_0$ in two occasions. A normalized SRD score close to 0 implies that the solution is similar to the reference---it is improbable that the resemblance happened by chance. On the other hand if the score is almost 1 then the solution vector ranks the components in a reverse order compared to the reference, which is equally unlikely.

Note that the SRD distribution depends heavily on the size of the input $n$ and the number of ties occurring in the solutions and the reference. Also, there is more than one meaningful way to derive the SRD distribution.

\cite{Diaconis_Graham_1977} provided a characterization for the SRD distribution when ties are not present. Let $S_n$ be the set of permutations of the first $n$ natural numbers and let $p, r \in S_n$ be chosen independently and uniformly. Furthermore let $D(p,r)$ be the random variable corresponding to the distance between $p$ and $r$ in $L_1$-norm. Then the following is true.

\begin{theorem}{(\cite{Diaconis_Graham_1977})\\}
If $n$ tends to infinity

\begin{gather}
    \mathbb{E}\left(D(p,r)\right)=\frac{1}{3}n^2+\mathcal{O}(n) \notag\\
    var\left(D(p,r)\right)=\frac{2}{45}n^3+\mathcal{O}(n^2) \notag
\end{gather}
\noindent and $D(p,r)$ follows asymptotically normal distribution.
\end{theorem}

This implies that when there are no ties, the expected value of a normalized SRD score of two random rankings is around $\frac{2}{3}$ (cf.\  Eq.~\ref{eq:normalization}). 

There are a few caveats here. The theorem does not help us when $n$ is small. In practice, the approximation error becomes sufficiently small when $n>13$. For smaller $n$ the exact distribution can be calculated by considering every permutation of length $n$. However, the presence of ties spoils the normality property and produces a deformed (zigzagged) probability density function.

In addition, it is not clear how the ties should be handled in the empirical SRD distribution. Should we fix the reference ranking and only generate random rankings representing the solutions? Should the number of ties in the random ranking follow the frequency of ties in the reference or the frequency of ties in the solution vectors? If there is one solution that has many ties, while the others do not, should this affect the CRRN test of all solutions or just the one which exhibits many ties? There is no universally correct answer to these questions. The rSRD package is specifically designed to accommodate such situations, allowing the practitioner to have control over how to handle them. It provides the flexibility for the practitioner to choose a distribution that best describes their specific use case, enabling them to make informed decisions based on their requirements and preferences.

\subsubsection{CVST}

Solutions with the same or very close SRD scores are not necessarily similar to each other. Rankings can differ from the reference in different sections and still be of the same distance. Cross-Validation combined with Statistical Testing (CVST) is designed to determine whether two solutions are inherently the same or not. Consequently, the null-hypothesis is formulated as follows. Let $r^i$ and $r^j$ be two rankings corresponding to solution $i$ and $j$, in other words $r^j$ is the $j$th column vector of $R$.\\

$H_0$: \emph{$r^i$ and $r^j$ come from the same distribution.}\\

Since it is impossible to draw a statistical conclusion by comparing two single SRD values, a sample is created using cross-validation. Randomly selected rows are discarded, and the SRD values are re-calculated on the remaining table. The SRD scores vary a bit depending on whether the omitted section agrees with the reference or not. In this way, we can assign a set of SRD values to the solutions and infer their relationship by comparing the samples. This is a delicate issue as we have to strike a balance between type 1 and type 2 errors. The literature suggests a sample size between 5 and 10 \citep{Hastie2009}. Increasing the sample size (\textit{i.e.}\ the number of folds), above 10 will increase the bias but lower the variance\footnote{Note that \cite{Hastie2009}'s results relate to a slightly different setting. In our framework, cross-validated folds the discarded rows might partially overlap, hence the validation process is closer to bootstrapping.}.  

The rSRD package features the three most widely used statistical tests that are coupled with cross-validation: Wilcoxon, Dietterich, and Alpaydin. The reason why three different test is provided is that they handle type 1 and type 2 situations drastically differently \citep{Sziklai2024}. Wilcoxon excels in type 2 scenarios but it is inefficient in type 1 situations. In other words, it effectively identifies when solutions come from different distributions, but sometimes differentiates between solutions that come from the same distribution. Dietterich and Alpaydin perform quite well in type 1 cases but fail badly in type 2 scenarios. Ultimately only the practitioner knows which type of error is more costly to her, hence we implemented all three tests in the package.

Describing the exact mathematical formulation of the three statistical tests as well as discussing their performance on various data structures goes well beyond the scope of the manuscript. A detailed description can be found in ref.\ \citep{Sziklai2024}.

\subsection{Demonstrative example}

In this section, we present a small case study to demonstrate SRD and its validation options. The data is courtesy of \href{https://eulytix.eu}{Eulytix}. For presentational reasons, figures were re-created and enhanced by Graphics Layout Engine \citep{GLE_2022}.

The major part of the EU's legislative work is done in the European Parliament's Committees. Members of the European Parliament (MEPs) draw up, amend, and adopt legislative proposals and work out reports to be presented to the plenary assembly. To boost the chance of acceptance, MEPs routinely approach their fellow members to co-sponsor their initiative. EP Committees suffer from partisanship much less than national assemblies. This allows MEPs to draw supporters from a wider pool, quite possibly from parties not belonging to the same ideological family. One way to identify potential allies is to analyze the language they are using. In the following example, the texts of the amendments of the Committee on Industry, Research and Energy (ITRE) are analyzed in the 2014-2019 legislative term to profile MEPs. Expressions are classified into 16 categories (climate, infrastructure, security, etc.) and counted how many times an MEP used an expression belonging to a certain category in the amendments she co-sponsored. This reveals which areas the MEP deems important. Our assumption is that MEPs that have a similar profile are the most probable allies.

\begin{table}[t!]
\scriptsize
\begin{tabular}{lllllllll|l}
\hline
                                                               & Botenga & Bompard & Ernst & Chahim & Pereira & Buzek & Groothuis & Kaljurand & Rego \\ \hline
\begin{tabular}[c]{@{}l@{}}Budget/\\ Costs\end{tabular}        & 12      & 20      & 17    & 0      & 7       & 17    & 8         & 1         & 13   \\
\begin{tabular}[c]{@{}l@{}}Climate/\\ Environment\end{tabular} & 137     & 293     & 232   & 43     & 44      & 136   & 82        & 0         & 174  \\
Economy                                                        & 48      & 79      & 46    & 14     & 14      & 35    & 54        & 13        & 20   \\
Energy                                                         & 121     & 244     & 169   & 48     & 68      & 204   & 153       & 0         & 145  \\
Enterprise                                                     & 79      & 71      & 79    & 12     & 82      & 80    & 164       & 8         & 82   \\
Future                                                         & 17      & 66      & 65    & 4      & 11      & 35    & 75        & 40        & 16   \\
Geography                                                      & 53      & 121     & 85    & 13     & 28      & 49    & 39        & 13        & 44  \\
Health                                                         & 135     & 168     & 39    & 1      & 52      & 0     & 6         & 0         & 107  \\
Industry                                                       & 99      & 211     & 105   & 25     & 47      & 75    & 63        & 4         & 107  \\
Infrastructure                                                 & 14      & 50      & 32    & 10     & 17      & 54    & 83        & 5         & 24   \\
Labour                                                         & 66      & 111     & 84    & 9      & 40      & 3     & 5         & 2         & 40   \\
\begin{tabular}[c]{@{}l@{}}Mining/\\ Resources\end{tabular}    & 27      & 155     & 138   & 11     & 0       & 112   & 50        & 0         & 115  \\
Mode of conduct                                                & 35      & 48      & 34    & 9      & 6       & 7     & 16        & 13        & 28   \\
Science                                                        & 47      & 59      & 28    & 6      & 21      & 13    & 30        & 0         & 46   \\
Security                                                       & 18      & 67      & 39    & 1      & 5       & 18    & 27        & 76        & 20   \\
Social                                                         & 65      & 74      & 80    & 16     & 37      & 16    & 13        & 4         & 49   \\ \hline
$\underline{SRD}_i$ & 0.234 & 0.297 & 0.312 & 0.352 & 0.352 & 0.484 & 0.547 & 0.891 & 0
\end{tabular}
\caption{Profiles of the Members of the European Parliament (MEP). Breakdown of the number of expressions used by MEPs in amendments they co-sponsored by categories. MEP Sira Rego's profile serves as a reference. The last row displays the derived SRD scores.}\label{tab:itre_input}
\end{table}

Table~\ref{tab:itre_input} compiles the language profile of 9 MEPs of the ITRE Committee. The last column shows the number of expressions the MEP Sira Rego used in the amendments she co-sponsored. The last row shows the normalized SRD scores. MEP Botenga's profile is the closest to the reference, which is MEP Rego's language-pattern. Both of them find Climate/Environment the most important while contribute to Budget/Costs category the least relatively to other topics. In general, they rank the topics the same way, which suggests that they might be ideologically close.

\begin{figure}[t]
\centering
\includegraphics[width=12cm]{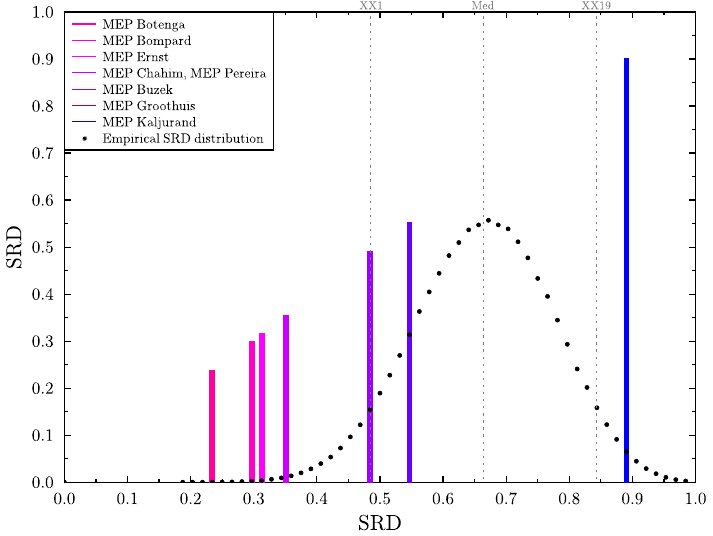}
\caption{CRRN test: The colored bars from left to right follow the same order as the legend from top to bottom, and their heights are equal to their normalized SRD values, hence expressing the distance from zero. The further they fall from the origin, the less they resemble to the reference. The vertical lines, XX1 and XX19 correspond to the 5\%  and 95\% threshold respectively.}\label{fig:itre_crrn}
\end{figure}

Figure~\ref{fig:itre_crrn} shows the result of the CRRN test. For sake of simplicity, SRD distribution was generated by assuming there are no ties, although this is not true. Both the reference and the solution vectors contain tied values. Six MEPs, Botenga, Bompard, Ernst, Chahim, Pereira, and Buzek pass the test, meaning the SRD value corresponding to their profile is below the 5\% significance threshold, marked by the first dashed line labeled XX1.

There are a few things worth noting. First, MEP Buzek's SRD value falls exactly on the 5\% threshold. Had we considered the ties when generating the SRD distribution, Buzek's profile would have fallen to the right of the threshold\footnote{Technically, this makes no difference as we cannot reject the null hypothesis on either case. However, it shows that with enough data, some bound to fall near to the threshold, and a change in the distribution could mean pushing a solution over the threshold.}. Second, MEP Kaljurand's profile is so different from Rego's, that it cannot even be considered random, but rather Rego's exact opposite. Reverse orderings are always revelatory, in our case it probably indicates that Kaljurand's and Rego's worldview is not compatible. This should be taken with a grain of salt as we have much fewer observations for Kaljurand, than for other MEPs. 
Groothuis' SRD score falls between the 5\% and 95\% significance thresholds meaning it cannot be distinguished from a random ranking. That suggests, that Groothuis is neither a probable ally nor an adversary of Rego.
Finally, MEPs Chahim and Pereira's rankings are of the same distance from the reference and Bompard's and Ernst's SRD scores are also very close. This calls for cross-validation.

\begin{figure}[t]
\centering
\includegraphics[width=12cm]{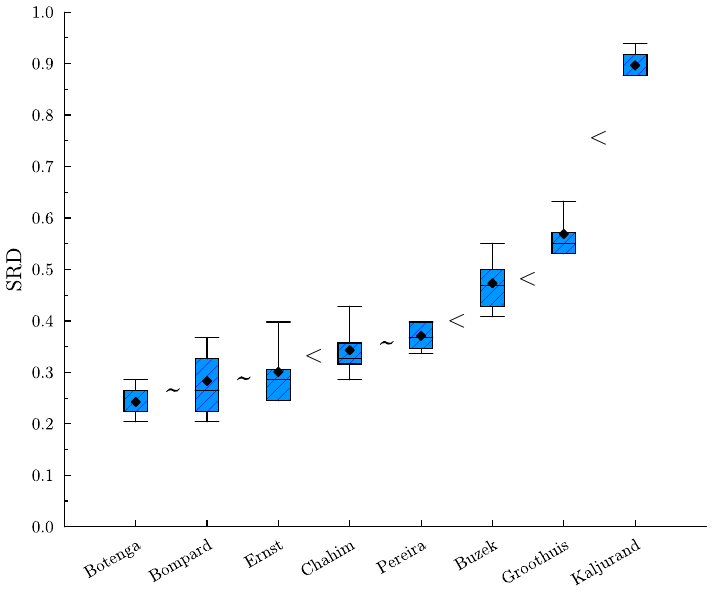}
\caption{Cross-validation: Box whisker plot representing the SRD values calculated on the different folds. The $<$ symbol denotes significant difference between the values, while $\sim$ marks that the null-hypothesis could not be rejected. }\label{fig:itre_boxplot}
\end{figure}

Figure~\ref{fig:itre_boxplot} displays the result of the cross-validation coupled with an 8-fold Wilcoxon test. We took 8 samples from the input matrix, by repeatedly removing two random rows, calculated the SRD values on all the samples, then ordered the solutions (\textit{i.e.}\ the MEPs) by average SRD, finally performed a Wilcoxon test between the subsequent solution-pairs. The test found that there is no significant difference between Bompard and Ernst and Chahim and Pereira. The choice of the test matters. Wilcoxon is much more sensitive than Alpaydin or the Dietterich test. A 10-fold Alpaydin-test finds only the last difference (between Groothuis and Kaljurand) significant, while a 10-fold Dietterich finds no significant difference between the \textit{consecutive} solution-pairs. 

\section{Package features}

In this section, we review the features of the package through examples. The functions of rSRD are divided into three categories. Core functions that relate to SRD computation and validation are given a name that starts with the '\texttt{calculate}' prefix. Plot generating functions start with the '\texttt{plot}', while utility functions start with the '\texttt{utils}' prefix.

Throughout this section, we will demonstrate package features using the MEP profiles we introduced in the previous section and a football dataset compiled from \url{https://www.whoscored.com/}. Table~\ref{tab:bundesliga} presents various game statistics aggregated over the 2020/21 season of the Bundesliga. The last column displays the points gathered by the teams. The final position in the Tableau is a natural reference that encompasses how strongly a team performed during the season. Some game statistics closely follow, while others seem to be independent of this ranking. With an SRD analysis we can uncover the style of play in the Bundesliga and explore which game elements dominate and which ones are not pertinent to the success of the teams. 
In the following, we will use '\texttt{mep\_profiles.csv}' and '\texttt{football\_leagues.csv}' to demonstrate package features.

\begin{verbatim}
R> profiles_df <- read.csv("mep_profiles.csv", row.names = 1,
   sep = ";")
R> bundesliga_df <- read.csv("bundesliga20_21.csv", row.names = 
   1, sep = ";")
\end{verbatim}

%
\begin{table}[t]
\scriptsize
\begin{tabular}{lccccccc|c}
\hline
Teams         & \multicolumn{1}{c}{Shots pg} & \multicolumn{1}{c}{RY cards} & \multicolumn{1}{c}{Possession\%} & \multicolumn{1}{c}{Pass} & \multicolumn{1}{c}{Dribbles pg} & \multicolumn{1}{c}{Offsides pg} & \multicolumn{1}{c|}{Fouls pg} & \multicolumn{1}{c}{pts} \\ \hline
Bayern        & \multirow{2}{*}{19.8}        & \multirow{2}{*}{38}          & \multirow{2}{*}{64.8}            & \multirow{2}{*}{86}      & \multirow{2}{*}{14.5}           & \multirow{2}{*}{2.2}            & \multirow{2}{*}{9}            & \multirow{2}{*}{77}     \\
Muenchen        &                              &                              &                                  &                          &                                 &                                 &                               &                         \\
Bayer         & \multirow{2}{*}{13.5}        & \multirow{2}{*}{66}          & \multirow{2}{*}{53.7}            & \multirow{2}{*}{81.8}    & \multirow{2}{*}{11.8}           & \multirow{2}{*}{1.9}            & \multirow{2}{*}{10.7}         & \multirow{2}{*}{64}     \\
Leverkusen    &                              &                              &                                  &                          &                                 &                                 &                               &                         \\
Borussia      & \multirow{2}{*}{13.3}        & \multirow{2}{*}{62}          & \multirow{2}{*}{59.4}            & \multirow{2}{*}{84}      & \multirow{2}{*}{10.4}           & \multirow{2}{*}{2.1}            & \multirow{2}{*}{10.6}         & \multirow{2}{*}{69}     \\
Dortmund      &                              &                              &                                  &                          &                                 &                                 &                               &                         \\
RB Leipzig    & \multicolumn{1}{c}{12.9}     & \multicolumn{1}{c}{49}       & \multicolumn{1}{c}{56.5}         & \multicolumn{1}{c}{83.1} & \multicolumn{1}{c}{10.2}        & \multicolumn{1}{c}{1.9}         & \multicolumn{1}{c|}{10.6}     & \multicolumn{1}{c}{58}  \\
SC Freiburg      & \multicolumn{1}{c}{13.6}     & \multicolumn{1}{c}{34}       & \multicolumn{1}{c}{48.6}         & \multicolumn{1}{c}{76.2} & \multicolumn{1}{c}{6.6}         & \multicolumn{1}{c}{1.7}         & \multicolumn{1}{c|}{11.5}     & \multicolumn{1}{c}{55}  \\
Borussia      & \multirow{2}{*}{14.8}        & \multirow{2}{*}{67}          & \multirow{2}{*}{54.1}            & \multirow{2}{*}{82}      & \multirow{2}{*}{10.9}           & \multirow{2}{*}{2.2}            & \multirow{2}{*}{10.6}         & \multirow{2}{*}{45}     \\
M.Gladbach    &                              &                              &                                  &                          &                                 &                                 &                               &                         \\
1.~FC Koeln       & \multicolumn{1}{c}{13.8}     & \multicolumn{1}{c}{67}       & \multicolumn{1}{c}{54.8}         & \multicolumn{1}{c}{77.4} & \multicolumn{1}{c}{7.4}         & \multicolumn{1}{c}{2}           & \multicolumn{1}{c|}{12.3}     & \multicolumn{1}{c}{52}  \\
FSV Mainz 05       & \multicolumn{1}{c}{13.8}     & \multicolumn{1}{c}{62}       & \multicolumn{1}{c}{46}           & \multicolumn{1}{c}{74.1} & \multicolumn{1}{c}{7.9}         & \multicolumn{1}{c}{1.9}         & \multicolumn{1}{c|}{14.6}     & \multicolumn{1}{c}{46}  \\
VfL Wolfsburg     & \multicolumn{1}{c}{12.4}     & \multicolumn{1}{c}{61}       & \multicolumn{1}{c}{50.2}         & \multicolumn{1}{c}{78.6} & \multicolumn{1}{c}{9.4}         & \multicolumn{1}{c}{2.2}         & \multicolumn{1}{c|}{12.1}     & \multicolumn{1}{c}{42}  \\
VfB Stuttgart & \multicolumn{1}{c}{13.3}     & \multicolumn{1}{c}{64}       & \multicolumn{1}{c}{50.4}         & \multicolumn{1}{c}{80.7} & \multicolumn{1}{c}{11}          & \multicolumn{1}{c}{1.6}         & \multicolumn{1}{c|}{10.9}     & \multicolumn{1}{c}{33}  \\
Union         & \multirow{2}{*}{12.1}        & \multirow{2}{*}{62}          & \multirow{2}{*}{43.3}            & \multirow{2}{*}{73.6}    & \multirow{2}{*}{7.2}            & \multirow{2}{*}{1.9}            & \multirow{2}{*}{12.4}         & \multirow{2}{*}{57}     \\
Berlin        &                              &                              &                                  &                          &                                 &                                 &                               &                         \\
Eintracht     & \multirow{2}{*}{13.2}        & \multirow{2}{*}{60}          & \multirow{2}{*}{49.4}            & \multirow{2}{*}{76.2}    & \multirow{2}{*}{8.4}            & \multirow{2}{*}{1.6}            & \multirow{2}{*}{12.4}         & \multirow{2}{*}{42}     \\
Frankfurt     &                              &                              &                                  &                          &                                 &                                 &                               &                         \\
TSG Hoffenheim    & \multicolumn{1}{c}{13.3}     & \multicolumn{1}{c}{75}       & \multicolumn{1}{c}{53.2}         & \multicolumn{1}{c}{80.7} & \multicolumn{1}{c}{7.4}         & \multicolumn{1}{c}{2.1}         & \multicolumn{1}{c|}{12.6}     & \multicolumn{1}{c}{46}  \\
VfL Bochum        & \multicolumn{1}{c}{12.1}     & \multicolumn{1}{c}{55}       & \multicolumn{1}{c}{44.5}         & \multicolumn{1}{c}{72.1} & \multicolumn{1}{c}{6.8}         & \multicolumn{1}{c}{2.2}         & \multicolumn{1}{c|}{12.5}     & \multicolumn{1}{c}{42}  \\
FC Augsburg      & \multicolumn{1}{c}{10.8}     & \multicolumn{1}{c}{74}       & \multicolumn{1}{c}{40.6}         & \multicolumn{1}{c}{72}   & \multicolumn{1}{c}{7.7}         & \multicolumn{1}{c}{2.1}         & \multicolumn{1}{c|}{13.2}     & \multicolumn{1}{c}{38}  \\
Arminia       & \multirow{2}{*}{10.7}        & \multirow{2}{*}{55}          & \multirow{2}{*}{39.9}            & \multirow{2}{*}{71.7}    & \multirow{2}{*}{8.4}            & \multirow{2}{*}{1.6}            & \multirow{2}{*}{12.7}         & \multirow{2}{*}{28}     \\
Bielefeld     &                              &                              &                                  &                          &                                 &                                 &                               &                         \\
Hertha BSC       & \multirow{2}{*}{10.8}        & \multirow{2}{*}{64}          & \multirow{2}{*}{43.2}            & \multirow{2}{*}{74.7}    & \multirow{2}{*}{8.3}            & \multirow{2}{*}{2.1}            & \multirow{2}{*}{12.4}         & \multirow{2}{*}{33}     \\
Berlin        &                              &                              &                                  &                          &                                 &                                 &                               &                         \\
Greuther      & \multirow{2}{*}{9.2}         & \multirow{2}{*}{61}          & \multirow{2}{*}{43}              & \multirow{2}{*}{74.8}    & \multirow{2}{*}{7.8}            & \multirow{2}{*}{2}              & \multirow{2}{*}{12.9}         & \multirow{2}{*}{18}     \\
Fuerth        &                              &                              &                                  &                          &                                 &                                 &                               &                         \\ \hline
\end{tabular}

\caption{Game statistics of Bundesliga teams aggregated over the season 2020/21 (pg stands for 'per game', RY abbreviates red and yellow).}\label{tab:bundesliga}
\end{table}

\subsection{Preprocessing}

The rSRD package offers a variety of utility functions. The Sum of Ranking Differences relies on the existence of a reference vector. One issue that often comes up in applications is the need for creating a reference. If an external reference (e.g.\ industry standard, known theoretical value) is not readily available, we might still be able to extract one from the data. For instance, if the solutions represent independent measurements then taking the average of the columns (\textit{i.e.}\ the average of the row values for each row) may work as a reference as the random errors cancel out. The median is often preferred in the presence of outliers. If objects represent performance properties that the solution needs to optimize, then row minimum or maximum might be a good choice. Yet sometimes we have a mix of these: some rows correspond to measurements, while other rows are performance traits.

In such cases, there are two issues we need to resolve. The first issue is that objects are not necessarily of the same kind; hence, their values do not fall on the same scale. For instance, the source data for the football leagues dataset measured the number of red and yellow cards in bulk for the whole season, while the number of offsides was averaged over the number of games. Even if the objects are on the same scale, it might be beneficial to standardize the values. For instance, in the MEP profiles, not every representative was equally productive. MEPs who submitted more amendments inevitably generated more phrases in each category, thus the absolute numbers are a bit misleading.

For the sake of convenience, rSRD includes standardization functions. Four methods are implemented: \texttt{scale\_to\_unit, standardize, range\_scale} and \texttt{scale\_to\_max} to meet the different needs.

\begin{description}
    \item [\texttt{scale\_to\_unit}] transforms each column vector into a vector of length 1.
    \item [\texttt{standardize}] converts the values of each column vector into a standard normal scale.
    \item [\texttt{range\_scale}] casts the values of each column vector into the $[0,1]$ interval.
    \item [\texttt{scale\_to\_max}] divides each column vector by the maximum value of that column.
\end{description}

For instance, the practitioner may prefer to use \texttt{standardize} if there are both positive and negative values present in the data, or use \texttt{range\_scale} if she wants to enforce each column into the same interval. Note that \texttt{utilsPreprocessDF()} transforms the columns. If the user would like to standardize the rows, the data frame should be transposed first.

\begin{verbatim}
R> utilsPreprocessDF(profiles_df, method="range_scale")
\end{verbatim}

After range scaling the MEP profiles, it is immediately clear what topic is the most important for each representative. In the case of MEP Chahim, it is Energy, while the value $0.25$ under Enterprise indicates that Chahim expressed a quarter of the amount of energy-related phrases about enterprises.

It is important to note, that each data preprocessing method is strictly monotonic. Hence, if we already have a fixed reference, the ranking transformation of SRD makes preprocessing \emph{by columns} redundant, meaning we will obtain the same SRD values with or without using preprocessing. However, if we normalize the rows, or we obtain the reference values as a function of the columns then normalization alters the outcome, and different normalization processes may lead to different SRD scores.

The second issue is the extraction of a reference when the rows correspond to different properties. The rSRD package makes reference creation a flexible process. The \texttt{utilsCreateReference()} method appends a new column at the end of the data frame. The \texttt{method} parameter describes the aggregation process. Five possibilities are available: \texttt{max, min, median, mean} and \texttt{mixed}. In the case of the first four, the row maximum, row minimum, row median, or row mean value is appended respectively.

\begin{verbatim}
R> SRD_input <- data.frame(
  A=c(2, 5, 7, 8),
  B=c(5, 1, 6, 10),
  C=c(6, 3, 2, 3))
\end{verbatim}

\begin{verbatim}
R> utilsCreateReference(SRD_input, method = "mean")
\end{verbatim}
\begin{verbatim}
  A  B C  refCol
1 2  5 6   4.33
2 5  1 3   3.00
3 7  6 2   5.00
4 8 10 3   7.00
\end{verbatim}

If the method is set to \texttt{mixed}, then we can specify an aggregation process for each row separately. For instance,

\begin{verbatim}
R> ref <- c("max","min","mean","mean")
R> utilsCreateReference(SRD_input, method = "mixed", ref)
\end{verbatim}

\begin{verbatim}
  A  B C refCol
1 2  5 6    6
2 5  1 3    1
3 7  6 2    5
4 8 10 3    7
\end{verbatim}

creates a reference for \texttt{SRD\_input} by taking the maximum of the first row, the minimum of the second, and the mean of the last two. The size of the \texttt{ref} vector must match the size of the input. Note that the function does not change the input, only returns with a new matrix.

\smallskip

rSRD provides a function for detailed SRD calculation, displaying the ranking transformation, the distance calculation, and the raw (unnormalized) SRD scores.

\begin{verbatim}
R> SRD_input = utilsCreateReference(SRD_input, method = "mixed", ref)
R> utilsDetailedSRD(SRD_input)
\end{verbatim}
\begin{verbatim}
  A  A_Rank A_Dist  B B_Rank B_Dist  C  C_Rank C_Dist  refCol refCol_Rank
1 2     1     2     5      2     1   6      4    1.0      6        3
2 5     2     1     1      1     0   3    2.5    1.5      1        1
3 7     3     1     6      3     1   2      1    1.0      5        2
4 8     4     0    10      4     0   3    2.5    1.5      7        4
5 -     -     4     -      -     2   -      -    5.0      -        -
\end{verbatim}

This utility serves multiple purposes---it allows users to visually check which part of the solution differs from the reference. Additionally, it is beneficial for pedagogical reasons, enabling users to quickly comprehend how the statistics are computed. Since the output of this function displays unnormalized SRD values, the user can promptly query the normalizing factor by calling

\begin{verbatim}
R> utilsMaxSRD(4)
\end{verbatim}
\begin{verbatim}
8
\end{verbatim}

which just returns the result of  Eq. \ref{eq:normalization}.
The rSRD package also allows to perform the ranking transformation without the SRD calculation. This option was built in to ensure compatibility with other ranking-based packages. In this way, the user can export the ranking matrix without the need to format the output.

\begin{verbatim}
R> utilsRankingMatrix(SRD_input)
\end{verbatim}
\begin{verbatim}
  A B   C
1 1 2  4.0
2 2 1  2.5
3 3 3  1.0
4 4 4  2.5
\end{verbatim}

\subsection{Core functions}

There are three core functions, one corresponding to the SRD calculation and the other two for the validation steps. Each function is capable of exporting the results to a CSV file.
First, let us look at the SRD values of the Bundesliga data.

\begin{verbatim}
R> calculateSRDValues(bundesliga_df)
\end{verbatim}
\begin{verbatim}
[1] 0.3395062 0.7037037 0.3148148 0.3950617 0.6049383 0.6604938 
    0.8888889
\end{verbatim}

The output follows the order of the columns in the input, in this case the column order of Table~\ref{tab:bundesliga}. The scores range over almost the entire $[0,1]$ interval. 
The number of shots is a good indicator of the final position in the tableau, as it ranks the football teams similarly as the obtained points. On the other hand, fouls committed per game seem to rank the teams in reverse order: worse teams tend to have more fouls. Which one of these scores are significant?

\begin{verbatim}
R> calculateSRDDistribution(bundesliga_df, option = "f")
\end{verbatim}
\begin{verbatim}
SRD_Distribution
    SRD_value   relative_frequency
1      0.0000        0.000000
2      0.2099        0.000001
3      0.2222        0.000001
4      0.2284        0.000005
...
127    0.9877        0.000023
128    0.9938        0.000003
129    1.0000        0.000019

xx1
[1] 0.4938
q1
[1] 0.5926
median
[1] 0.6667
q3
[1] 0.7346
xx19
[1] 0.8272
avg
[1] 0.6631711
std_dev
[1] 0.1020909
\end{verbatim}

The XX1 and XX19 values specify the 5\% and 95\% significance thresholds. Thus, the number of shots, ball possession, successful passes, and fouls committed are statistically significant. The first three rank the teams similarly to the reference $(\underline{SRD}_i<XX1)$, while fouls rank teams in reverse order $(\underline{SRD}_i\ge XX19)$. The function \texttt{calculateSRDDistribution} has four parameters, of which only the first, the input dataframe, is compulsory. The second parameter specifies how the SRD distribution should be generated. The available options are the following:

\begin{description}
    \item['n'] There are no ties for the solution vectors, the reference vector is fixed.

    \item['r'] There are no ties. Both the column vector and the reference are generated randomly.

    \item['t'] Ties occur with a fixed probability specified by the user for both the solution vectors and the reference vector.

    \item['p'] Ties occur with a fixed probability specified by the user for the solution vectors, the reference vector is fixed.

    \item['d'] Tie distribution reflects the tie frequencies displayed by the solution vectors, the reference vector is fixed.

    \item['f' (default)] Tie distribution reflects the tie frequencies displayed in the reference, the reference vector is fixed.
\end{description}

In each case, one million data points are produced by generating (either randomly or deterministically) solution and reference vectors and calculating their corresponding SRD scores. The third parameter, '\texttt{tie\_probability}' only plays a role if either options 't'  or 'p' were chosen, in which case they specify the tie frequencies occurring in the generated vectors. The package offers a utility function to check the number of ties present in a vector.

\begin{verbatim}
R> solution <- c(1,3,3,3,2,2,4,3)
R> utilsTieProbability(solution)
\end{verbatim}
\begin{verbatim}
[1] 0.5714286
\end{verbatim}

Note that in an $n$-long vector, ties can occur in $n-1$ places. In the above example, the values are ordered first to create a ranking. Hence, 4 out of the 7 places were tied and the result is 4/7 = 0.5714286.

If option 'r' is chosen, then both the solution and the reference are generated randomly without ties. In this case, the empirical distribution follows the Spearman's footrule distribution.

The choice of how to generate the distribution can have a huge impact on its probability density function. Although the values of XX1 and XX19 are more robust, small changes can mean that solutions falling near the significance threshold can gain or lose significance based on the selected method of generation\footnote{Although many authors warn that placing too much emphasis on the threshold is unreasonable and that one should simply report the $p$-value \citep{Wasserstein2016,Berrar2022}, the economic literature continues the practice of only accepting an effect when the corresponding variable is statistically significant.}. Therefore, the prudent way to design a hypothesis test is to select the generating method \textit{before} calculating the SRD values. Consequently, it is worthwhile to contemplate which option is suitable for which circumstances.

If the reference is based on observations like in the case of MEP profiles or the points obtained by football teams during a season, then its realization is prone to change with repeated measurements. The number of observed ties would vary with each measurement, thus there is no reason to treat the reference during the distribution generation fixed (see options 'r' or 't'). If, however, the reference is derived from established behavior, like for instance from theoretical properties of a substance, then the distribution should reflect the fact that the reference will not be changing (options 'n', 'p', 'd' or 'f').

If the compared solutions come from the same distribution or their distributions are highly similar, we are not committing a significant error by generating the solutions in the same way. However, if the solutions display highly different tie frequencies, we introduce a slight bias when comparing them under the same distribution. The question is how we prefer to represent the underlying solution space. Suppose there are two solutions with tie frequencies $x$ and $y$, where $x<<y$. Then we may consider generating a solution using the fixed tie probability $\frac{x+y}{2}$ (options 't' or 'p'). However, if the distance between $x$ and $y$ is too large, there is a chance that none of the generated vectors display a tie frequency as low as $x$ and as large as $y$. Another option would be to generate half the solutions with a tie frequency $x$ and the other half with a tie frequency $y$ (option 'd'), although in such cases we will end up with a hybrid distribution. If we have reason to believe that the solutions approximate or converge to the reference, using the same tie frequency for the solutions as for the reference is an adequate approach (default option). In such way, we will have an accurate estimation of the significance of solutions that have very low or very high SRD values. The point is, that we are not really interested in the solutions that behave randomly compared to the reference and exhibit an SRD value close to the expected value of the distribution.

Note that a fixed tie frequency does not mean that each vector is generated with the same number of ties. It simply means that each consecutive element in the vector is tied with the given probability.
Finally, let us stress that these are recommendations and ultimately only the practitioner can decide which distribution fits best to her use case.

To take a closer look at the cross-validation function of the package, let us consider again the Bundesliga dataframe. The number of shots, ball possession, and the ratio of successful passes are all statistically significant elements with relatively close SRD values. Ball possession seems to be the most significant factor, followed by the number of shots and successful passes. Can we be definitive about their order? Let us compare the solutions based on the results of the Wilcoxon and Alpaydin tests.

\begin{verbatim}
R> cv_Wilc <- calculateCrossValidation(bundesliga_df)
R> plotCrossValidation(cv_Wilc)
R> cv_Alp <- calculateCrossValidation(bundesliga_df, method = 
   "Alpaydin", number_of_folds = 10)
R> plotCrossValidation(cv_Alp)
\end{verbatim}

\begin{figure}[H]
\centering
\includegraphics[width=14cm]{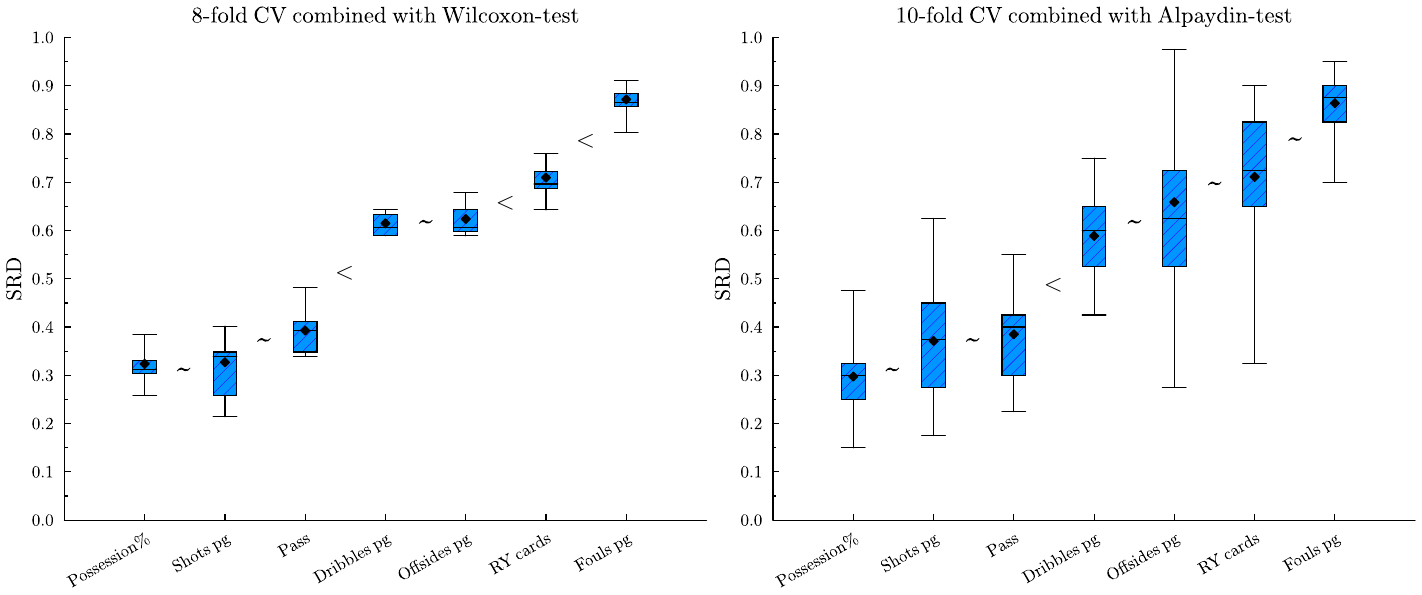}
\caption{CVST. Comparing the 8-fold Wilcoxon test with the 10-fold Alpaydin test. }\label{fig:bliga_boxplot}
\end{figure}

By default, the \texttt{calculateCrossValidation} function employs an 8-fold cross-validation combined with the Wilcoxon-test. Using the \texttt{method} and the \texttt{number\_of\_folds} parameters we can change the type of test and the number of folds used in the cross-validation. Figure~\ref{fig:bliga_boxplot} displays the box-whiskers plot created from the SRD scores of the different folds. Note that the SRD values vary much more under the Alpaydin-test than under the Wilcoxon-test. The reason stems from how the test performs cross-validation. Under Wilcoxon with $k$-fold cross-validation $\lceil\frac{n}{k}\rceil$ rows are removed in each fold, while under Alpaydin (and under Dietterich too) half the rows are removed in each fold.

There are a couple of things worth noting. First, the result of the cross-validation is stochastic. Each time the function is run, new random rows are selected for each fold. To ensure reproducibility, the function automatically saves not only the results but also the whole computation, including the selected rows. This feature can be disabled by setting the \texttt{output\_to\_file} parameter to false. Some information is only saved to the CSV file, the function does not return the number and indices of the removed rows. Secondly, the function orders solution based on the median of the computed SRD scores. Thus, RY cards, which was the second column of our dataframe, was relayed as the last but one column, while the number of shots switched places with ball possession. Only consecutive column pairs are tested. Thus, ball
possession does not significantly differ from the number of shots, and the number of shot does not significantly differ from successful passes, but the relationship between ball possession and successful passes is not tested. Finally, significance is classified into three categories: \texttt{n.s.}\ (not significant), $(p<0.1)$ and $(p<0.05*)$. Smaller $p$ values are usually not meaningful when the number of rows is not too large, hence they are omitted.

\begin{verbatim}
R> cv_Wilc
\end{verbatim}
\begin{verbatim}
new_column_order_based_on_folds
[1] 3 1 4 5 6 2 7

test_statistics
[1] 4 29 36 6 34 36

statistical_significance
[1] "n.s."   "(p<0.1)"   "(p<0.05*)"   "n.s."   "(p<0.05*)"   
    "(p<0.05*)"
			
(the rest of the report can be found in the Appendix)
\end{verbatim}

\subsection{Plotting}

\begin{figure}
\centering
\includegraphics[width=12cm]{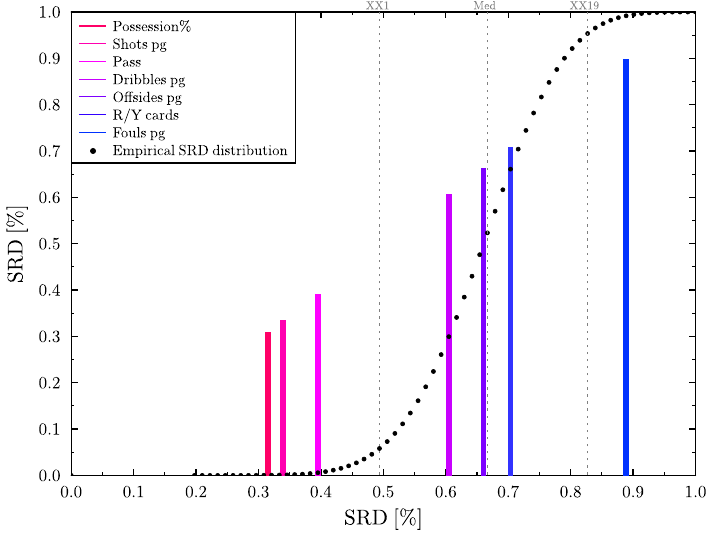}
\caption{CRRN test: SRD values are compared to the cumulative distribution of random rankings generated by method 'r'.}\label{fig:bliga}
\centering
\includegraphics[width=12cm]{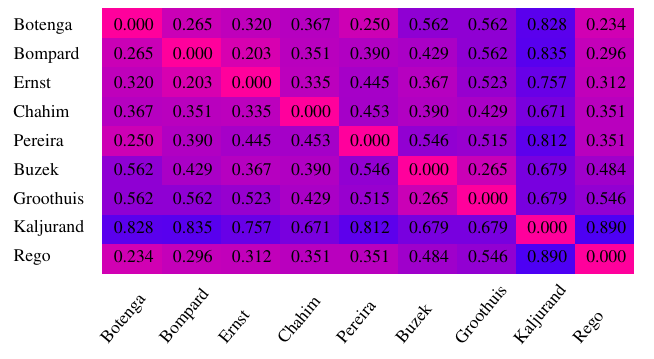}
\caption{Pairwise SRD distances.}\label{fig:heatmap}
\end{figure}

The package allows plotting the results of the permutation test under the chosen SRD distribution. In the previous section, we discussed the available distributions and their advantages. Here, we simply note that the package enables the user to choose between plotting the probability density function and the cumulative distribution function. If the \texttt{densityToDistr} parameter is set to true, then the SRD values are compared to the cdf (cf. \ref{fig:bliga}). Naturally, this does not affect the thresholds in any way.

\begin{verbatim}
R> dist_r <- calculateSRDDistribution(bundesliga_df, option = 'r')
R> plotPermTest(bundesliga_df, dist_r, densityToDistr = TRUE)
\end{verbatim}

In a multiple comparison setting, each solution might serve as a potential reference. In such cases, we are often interested in seeing how far the solutions deviate from each other. For instance, in the MEP profiles dataset, we might be interested in generating recommendations not only for MEP Sira Rego but for every representative. Setting each solution as the reference manually is a cumbersome task; thus, the package offers a way to do it in one step. The \texttt{plotHeatmapSRD} function not only calculates the pairwise distances and exports them to a CSV file, but also illustrates the distance matrix using a color palette. The package offers a built-in palette where a shade of red corresponds to an SRD score of 0, while a shade of blue corresponds to an SRD score of 1.

\begin{verbatim}
R> plotHeatmapSRD(profiles_df, output_to_file = TRUE, color = 
   utilsColorPalette)
\end{verbatim}

The color palette can be easily customized. The size of the palette indicates how many categories the $[0,1]$ interval is divided. For instance, the following code changes the color range from orange to green.

\begin{verbatim}
R> myPalette <- c("#eb9c34", "#ebba34", "#ebd634", "#ebe534", 
                  "#d9eb34", "#b7eb34", "#99eb34", "#6beb34")
R> plotHeatmapSRD(profiles_df, color = myPalette)
\end{verbatim}

\begin{figure}[t]
\centering
\includegraphics[width=10cm]{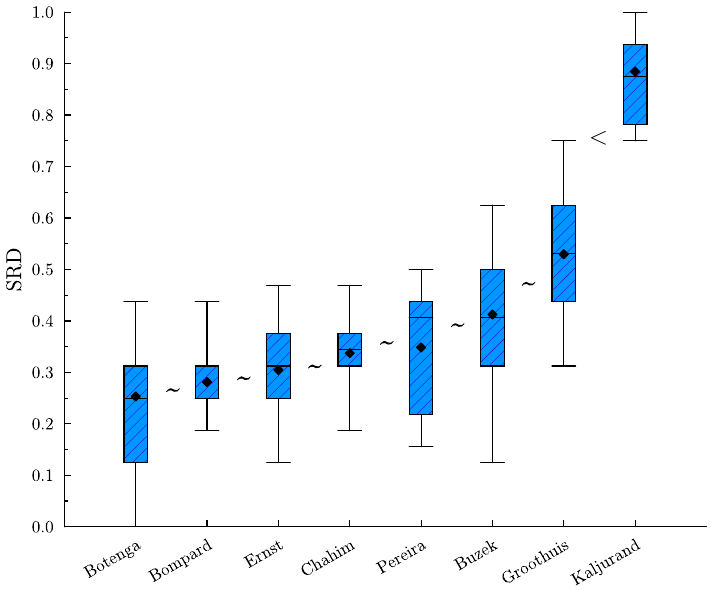}
\caption{CVST. Dietterich test on the MEP profiles. The $<$ symbol denotes significant difference between the values, while $\sim$ marks that the null-hypothesis could not be rejected. }\label{fig:mep_boxplot_d}
\end{figure}

Finally, the \texttt{plotCrossValidation} function produces a standard box-whisker plot (see Figure~\ref{fig:mep_boxplot_d}), where the whiskers mark the minimum and maximum of the SRD values calculated on the different folds. The box represents the first and third quartiles, the horizontal line inside the box is the median, while the crossmark/diamond symbol indicates the average.

\begin{verbatim}
cv <- calculateCrossValidation(profiles_df, method = "Dietterich", 
      number_of_folds = 10)
plotCrossValidation(cv)
\end{verbatim}

\section{Summary and discussion} \label{sec:summary}

Comparing solutions through a reference is a common task that appears in various fields of science, from marketing through sports to machine learning and material science. Sum of Ranking Differences is a conceptually simple, non-parametric statistical procedure that, besides ranking the solutions, can also separate the statistically significant ones. The popularity of this comparative analytical took stems from the fact that it was designed especially for comparing in the presence of a reference.

The rSRD package offers a comprehensive toolkit that encompasses data preprocessing, SRD computation, validation, and plotting features. Among the reviewed functions, validation is by far the most complex and demands the most consideration from the practitioner. The package features a variety of ways to compute the SRD distribution for the permutation test, and thus, can aid the practitioner in designing their hypothesis test. The choice may affect significance thresholds, so careful consideration is needed.

All validation options can export the calculations, which allows reproducibility. Apart from simplifying the validation process, the package also offers plotting features. In the future, the package will be updated based on user feedback. In particular, radar plots are planned to highlight how the SRD scores vary between different datasets.


Further updates aim to include an option to adjust the default sample size used for generating random rankings during the empirical distribution process. Currently, one million data points are generated, yielding precise results for all the examples tested so far. However, for very large $n$, the empirical distribution might be slightly inaccurate. Conversely, for smaller cases, generating one million data points is excessive, and reducing the sample size could speed up computation without sacrificing accuracy.

Another potential update pertains to cross-validation. Presently, cross-validation folds are generated by randomly discarding rows from the input table. In time-series analysis, it is common practice to exclude consecutive blocks rather than random rows. Incorporating block cross-validation, including predefined row sets, could be a straightforward and beneficial improvement.

Finally, the stochastic elements in the computation are implemented in \textsf{C++}. A drawback of this approach is that seeding the random number generator in \textsf{R} does not affect the computation. Although the calculations are saved into a CSV file, making the tests reproducible, incorporating a seed parameter into the cross-validation function would make the package more convenient.

\section*{Acknowledgement}
This work was supported by the Ministry of Innovation and Technology of Hungary from the National Research, Development and Innovation Fund, financed under the K type funding scheme, K 134260 (Héberger), K 146320 (Sziklai) and the FK funding scheme 137577 (Gere) . The last author thanks the Bavarian State Ministry of Science and Arts for partial funding.
\bibliographystyle{chicago}
\bibliography{bib}



\appendix
\section*{Appendix}

\begin{verbatim}
R> cv_Wilc

new_column_order_based_on_folds
[1] 3 1 4 5 6 2 7

test_statistics
[1] 4 29 36 6 34 36

statistical_significance
[1] "n.s."   "(p<0.1)"   "(p<0.05*)"   "n.s."   "(p<0.05*)"   
    "(p<0.05*)"

SRD_values_of_different_folds
	      Possession	 Shots	    Pass	   Dribbles	 Offsides
fold_1	0.3660714	0.3392857	0.4107143	0.6071429	0.6250000
fold_2	0.3839286	0.2589286	0.4821429	0.6428571	0.6428571
fold_3	0.3303571	0.3214286	0.3928571	0.6339286	0.6517857
fold_4	0.3035714	0.3482143	0.4107143	0.6428571	0.5982143
fold_5	0.3214286	0.3928571	0.3660714	0.5892857	0.5982143
fold_6	0.2589286	0.2142857	0.3392857	0.5892857	0.6785714
fold_7	0.3125000	0.3392857	0.3928571	0.6250000	0.5892857
fold_8	0.3125000	0.4017857	0.3482143	0.5892857	0.6071429

RY.cards	Fouls				
0.6964286	0.8660714				
0.6875000	0.8035714				
0.7589286	0.9107143				
0.7232143	0.8571429				
0.6875000	0.8839286				
0.6428571	0.9107143				
0.7589286	0.8660714				
0.7232143	0.8750000				




boxplot_values
	   Possession	Shots	Pass	Dribbles	Offsides
min	    0.2589	0.2143	0.3393	0.5893	0.5893
xx1	    0.2589	0.2143	0.3393	0.5893	0.5893
q1	     0.3036	0.2589	0.3482	0.5893	0.5982
median	 0.3125	0.3393	0.3929	0.6071	0.6071
q3	     0.3304	0.3482	0.4107	0.6339	0.6429
xx19	   0.3839	0.4018	0.4821	0.6429	0.6786
max	    0.3839	0.4018	0.4821	0.6429	0.6786

RY.cards	Fouls				
0.6429	0.8036				
0.6429	0.8036				
0.6875	0.8571				
0.6964	0.8661				
0.7232	0.8839				
0.7589	0.9107				
0.7589	0.9107				

\end{verbatim}

\end{document}